\documentclass[a4paper]{amsart}
\usepackage{graphicx}
\usepackage{latexsym}
\usepackage{hyperref}
\usepackage{amsfonts}
\usepackage[all]{xy}
\usepackage{amssymb, mathrsfs, amsfonts, amsmath}
\usepackage{amsfonts}
\usepackage{graphicx}
\usepackage{amscd}
\usepackage{amsmath}

\usepackage{enumerate}
\usepackage{graphicx}
\usepackage{eulervm}
\usepackage{geometry}
\usepackage{textcomp}
\usepackage{graphicx,color}

\newtheorem{theorem}{Theorem}

\newtheorem{corollary}{Corollary}

\newtheorem{definition}{Definition}

\newtheorem{lemma}{Lemma}

\newtheorem{proposition}{Proposition}
\newtheorem{remark}{Remark}

\newcommand{\R}{\mathbb{R}}

\newcommand{\diver}{{{\rm{div}\,}}}
\newcommand{\grad}{{{\rm{grad}\,}}}

\DeclareMathOperator{\Hess}{Hess}

\begin{document}
\title[On the Omori-Yau Maximum Principle and geometric applications]{On the Omori-Yau Maximum Principle and geometric applications}
{\author{Barnabe Pessoa Lima}
\address{Departamento de Matem\'{a}tica-UFPI\\
64049-550-Teresina-PI-Br} \email{barnabe@ufpi.edu.br} \thanks{The first author was partially supported by PROCAD-CAPES}}
{\author{Leandro de Freitas Pessoa}
\address{Departamento de Matem\'{a}tica-UFPI\\
64049-550-Teresina-PI-Br} \email{leandrofreitasp@yahoo.com.br} \thanks{The second author was partially supported by PICME-CAPES }}
\date{February 28, 2008}
\keywords{Omori-Yau maximum principle }
\footnote{{\it 2010 Mathematics Subject Classification:} Primary 53C42, ; Secondary 35B50}

\maketitle

\begin{abstract}We  introduce a version of the Omori-Yau maximum principle which generalizes the version obtained by Pigola-Rigoli-Setti \cite{prs-memoirs}. We apply our method to derive a non-trivial generalization  Jorge-Koutrofiotis Theorem  \cite{jorge-koutrofiotis} for cylindrically bounded submanifolds  due to Alias-Bessa-Montenegro \cite{alias-bessa-montenegro}, we extend results due to Alias-Dajczer \cite{alias-dajczer}, Alias-Bessa-Dajczer \cite{alias-bessa-dajczer} and Alias-Impera-Rigoli \cite{alias-impera-rigoli}.
\end{abstract}

\vspace{.2mm}

\section{Introduction and Statement of Results}

H. Omori \cite{omori}, studying isometric immersions of minimal submanifolds into cones of $\mathbb{R}^{n}$ proved the following global version of the maximum principle  for complete Riemannian manifolds with sectional curvature bounded below.

\begin{theorem}[Omori] \label{r1}  Let $M$ be a complete Riemannian manifold with sectional curvature  bounded below $K_{M}\geq -\Lambda^{2}$. If  $u\in C^{2}(M)$ with $u^{\ast}=\sup_{M}u<\infty$ then there exists a sequence of points $x_{k}\in M$, depending on $M$ and on $u$,  such that
\begin{eqnarray}\label{eq-omori}
\lim_{k\to \infty}u(x_{k}) =u^{\ast},& \vert \grad u\vert  (x_{k}) < \displaystyle \frac{1}{k},&  \Hess u(x_{k})(X,X) <  \displaystyle \frac{1}{k}\cdot \vert  X \vert^{2},
\end{eqnarray}for every $X\in T_{x_k}M$.
\end{theorem}

H. Omori's maximum principle was refined by S. T. Yau in a series of papers    \cite{yau}, \cite{yau2}, \cite{cheng-yau} (this later with S. Y. Cheng)  and applied to find elegant solutions to various analytic-geometric problems on Riemannian manifolds.  The version of the maximum principle Cheng-Yau proved  is the following variation of Theorem \ref{r1}.

\begin{theorem}[Cheng-Yau]\label{r2}
 Let $M$ be a complete Riemannian manifold with Ricci curvature  bounded below $Ric_{M}\geq -(n-1)\cdot\Lambda^{2}$. Then for any $u\in C^{2}(M)$ with $u^{\ast}=\sup_{M}u<\infty$, there exists a sequence of points $x_{k}\in M$, depending on $M$ and on $u$,  such that
\begin{eqnarray}\label{eq-omori}
\lim_{k\to \infty}u(x_{k}) =u^{\ast},& \vert \grad u\vert  (x_{k}) < \displaystyle \frac{1}{k},&  \triangle u(x_{k}) <  \displaystyle \frac{1}{k}.
\end{eqnarray}
\end{theorem}
H. Omori's  Theorem was  extended by C. Dias in \cite{dias} and Cheng-Yau's Theorem was extended by  Chen-Xin in \cite{chen-xin}. Recently,
 S. Pigola, M. Rigoli and A. Setti in their beautiful book \cite{prs-memoirs} introduced the following important concept.
 \begin{definition} The Omori-Yau maximum principle is said to hold on $M$ if for any given $u\in C^{2}(M)$ with $u^{\ast}=\sup_{M}u<\infty$, there exists a sequence of points $x_{k}\in M$, depending on $M$ and on $u$,  such that
\begin{eqnarray}\label{eq-omori-2}
\lim_{k\to \infty}u(x_{k}) =u^{\ast},& \vert \grad u\vert  (x_{k}) < \displaystyle \frac{1}{k},&  \triangle u(x_{k}) <  \displaystyle \frac{1}{k}.
\end{eqnarray} Likewise, the Omori-Yau maximum principle \textit{for the Hessian} is said to hold on $M$ if for any given $u\in C^{2}(M)$ with $u^{\ast}=\sup_{M}u<\infty$, there exists a sequence of points $x_{k}\in M$, depending on $M$ and on $u$,  such that \begin{eqnarray}\label{eq-omori}
\lim_{k\to \infty}u(x_{k}) =u^{\ast},& \vert \grad u\vert  (x_{k}) < \displaystyle \frac{1}{k},&  \Hess u(x_{k})(X,X) <  \displaystyle \frac{1}{k}\cdot \vert  X \vert^{2},
\end{eqnarray}for every $X\in T_{x_k}M$.
 \end{definition} This concept  was a new point of view of the Omori-Yau maximum principle. That is,  some geometries  do hold the Omori-Yau maximum principle whereas some does not.  That raised naturally the question of what are the geometries that hold the Omori-Yau maximum principle? The Omori-Yau maximum principle was shown to hold in several   geometric  settings, see for instance \cite{bps-revista}, \cite {fontenele-xavier}, \cite{lee-kim}, \cite{pigola-rigoli-rimoldi}, \cite{ratto-rigoli-setti}. Regarding this problem, S. Pigola, M. Rigoli and A. Setti \cite[pp. 7--10]{prs-memoirs} proved the a general class of Riemannian manifolds hold the  Omori-Yau maximum principle. They proved the following theorem.
\begin{theorem}[Pigola-Rigoli-Setti]\label{r3}
Let $M$ be a Riemannian manifold and assume that there exists a non-negative function $\gamma $ satisfying the following:
\begin{itemize}
\item[C1)] $\gamma (x) \rightarrow  +\infty $, \,\,\, as \,\,\, $x\rightarrow \infty $;
\item[]
\item[C2)] $\exists A > 0 \,\, \textit{such} \,\textit{that}  \,\, \vert \grad \gamma \vert  \leq  A \cdot \sqrt{\gamma} $, off a compact set;
\item[]
\item[C3)] $\exists B > 0 \,\, \textit{such} \,\textit{that} \,\, \triangle \gamma \leq B\cdot \sqrt{\gamma \,G(\sqrt{\gamma})} $, off a compact set;
    \item[] where $ G:[0,+\infty)\to [0, \infty) $ is a smooth function  satisfying
\item[] \begin{equation}\label{eq5}
\begin{array}{llll} G(0) > 0, &   G'(t) \geq 0,& \displaystyle{\int_{0}^{+\infty}\frac{ds}{\sqrt{G(s)}} = +\infty},&\displaystyle{\limsup_{t \rightarrow +\infty}\frac{tG(\sqrt{t})}{G(t)} < +\infty}.
\end{array}
\end{equation}
\end{itemize}
Then the Omori-Yau maximum principle holds on $M$.
\end{theorem} If instead of $C3)$ we assume the following stronger hypothesis
\begin{enumerate}
\item[C4)] $\exists B > 0 $ such that $ \Hess \gamma (.,.) \leq B\sqrt{\gamma G(\sqrt{\gamma})}\langle.,.\rangle $, off a compact set.
\end{enumerate}Then the Omori-Yau maximum principle for the Hessian holds on $M$.

\begin{remark}An example of a smooth function $ G\in C^{\infty}([0,+\infty)) $ satisfying (\ref{eq5}) is given by $G(t)=t^{2}\,\Pi_{j=1}^{\ell}(\log^{(j)}(t))^{2}$ for $t\gg 1$, where $\log^{(j)}$ is the $j$-th iterated logarithm and $\ell\in \mathbb{N}$.
\end{remark} \begin{remark}Quoting Pigola-Rigoli-Setti, \textquotedblleft The proof of Theorem \ref{r3} shows that one needs $\gamma$ to be $C^{2}$ only in a neighborhood of $x_{k}$. This is the case that $\gamma=\rho^{2}$ is the square of the Riemannian distance from a fixed point $x_{o}$ and $x_{k}$ is not on the cut locus of $x_{o}$. The case that $x_k$ is the cut locus of $x_o$ can be dealt with a trick of Calabi \cite{calabi} so that we may assume that $\gamma$ is always $C^{2}$ in a neighborhood of $x_k$.\textquotedblright,    \cite[Remark 1.11]{prs-memoirs}.
\end{remark}
\begin{corollary}[Pigola-Rigoli-Setti] Let $M$ be a complete Riemannian manifold with Ricci curvature satisfying $$ Ric (x)\geq - B^{2}\cdot G(\rho (x)), \, \, {\rm for}\,\, \rho (x)\gg 1$$ where $G\in C^{\infty}([0,+\infty)) $  satisfies  \eqref{eq5}, $\rho (x)={\rm dist}_{M}(x_{0}, x)$, $B\in \mathbb{R}$. Then $\gamma=\rho^{2}$ satisfies $C1-C3$. Therefore the Omori-Yau maximum principle  holds on $M$ by Theorem \ref{r3}. This shows that Theorem \ref{r3} extends Theorem \ref{r2}.
\end{corollary}

In this paper we give an extension of Pigola-Rigoli-Setti's Theorem \ref{r3}. We prove the following result.

\begin{theorem}[Main Theorem]\label{Principal Theorem}
Let $M$ be a complete Riemannian manifold and assume that there exists a non-negative function $\gamma $ satisfying:
\begin{itemize}
\item[h1)] $\gamma (x) \rightarrow  +\infty $, \,\,\, as \,\,\, $x\rightarrow \infty $;
\item[]
\item[h2)]$\exists A > 0 \,\, \textit{such}\,\textit{that}  \,\, \vert \grad \gamma\vert \leq A \cdot \sqrt{G(\gamma)}\left(\int_{0}^{\gamma}\frac{ds}{\sqrt{G(s)}}+1\right) $ off a compact set;
\item[]
\item[h3)]$\exists B > 0 \,\, \textit{such} \,\textit{that} \,\, \triangle \gamma \leq B \cdot \sqrt{G(\gamma)}\left(\int_{0}^{\gamma}\frac{ds}{\sqrt{G(s)}}+1\right) $ off a compact set;
     \item[] where $ G:[0,+\infty)\to [0, \infty) $ is a smooth function  satisfying
\item[] \begin{equation}\label{eq6}
\begin{array}{llll} G(0) > 0, &   G'(t) \geq 0,& \displaystyle{\int_{0}^{+\infty}\frac{ds}{\sqrt{G(s)}} = +\infty}.
\end{array}
\end{equation}
    \end{itemize} Then if $ u \in C^{2}(M)$
satisfies
$
\displaystyle{\lim_{x \rightarrow \infty}\frac{u(x)}{\varphi(\gamma(x))}} = 0 ,
$
where
$
\displaystyle{\varphi(t) = \log\left(\int_{0}^{t}\frac{ds}{\sqrt{G(s)}}+1\right)},
$
then there exists a sequence ${x_{k}} \in M$, $k \in \mathbb{N} $ such that
\begin{enumerate}
\item[]\begin{equation} \begin{array}{lll}\displaystyle{\vert \grad u\vert(x_{k}) < \frac{1}{k}},& &
 \displaystyle{\triangle u(x_{k}) < \frac{1}{k}} \end{array}\end{equation}
\end{enumerate}\end{theorem}If instead of $h3)$ we assume
$ \displaystyle{\Hess \gamma (.,.)  \leq \sqrt{G(\gamma)}\left(\int_{0}^{\gamma}\frac{ds}{\sqrt{G(s)}}+1\right)}\langle .,. \rangle $ off a compact set, then
 $ \displaystyle{\Hess u(x_{k}) ( .,. ) < \frac{1}{k} \langle . , . \rangle} $.
Moreover, if $u \in C^{2}(M)$  is bounded above $u^{\ast}=\sup_{M}u<\infty$ then  $u(x_{k})\to u^{\ast}$.

This result above should be compared with  a fairly recent result \cite[Cor. A1.]{prs-jfa} also due to S. Pigola, M. Rigoli, A. Setti improved Theorem \ref{r3} to more general elliptic operators.

\vspace{.2mm}

\section{Omori-Yau maximum principle}\label{sec-main}

\textbf{Proof of Theorem \ref{Principal Theorem}:}
We fix a sequence of positive real numbers $(\varepsilon_{k})_{k \in \mathbb{N}}$ such that, $ \varepsilon_{k} \rightarrow 0$ and \linebreak
consider now any function $u \in C^{2}(M)$ satisfying $
\displaystyle{\lim_{x \rightarrow \infty}\frac{u(x)}{\varphi(\gamma(x))}} = 0
$, where \linebreak
$
\displaystyle{\varphi(t) = \log\left(\int_{0}^{t}\frac{ds}{\sqrt{G(s)}}+1\right)}
$. Define
\begin{eqnarray}\label{defgk}
g_{k}(x) = u(x) - \varepsilon_{k}\varphi(\gamma(x))
\end{eqnarray}
and observe that $\varphi $ is $C^{2}(M)$, positive and satisfies
\begin{eqnarray*}
\varphi(t) \rightarrow +\infty  \,\,\,\, as \,\,\,\, t\rightarrow  +\infty .
\end{eqnarray*}
By a direct computation we have
\begin{eqnarray}\label{filinha}
\varphi'(t)&=& \left[\sqrt{G(t)}\left(\int_{0}^{t}\frac{ds}{\sqrt{G(s)}} + 1\right)\right]^{-1},\nonumber \\
&& \nonumber \\
\varphi''(t) &=& - \left[\sqrt{G(t)}\left(\int_{0}^{t}\frac{ds}{\sqrt{G(s)}} + 1\right)\right]^{-2}\left[\frac{G'(t)}{2\sqrt{(G(t))}}\left(\int_{0}^{t}\frac{ds}{\sqrt{G(s)}} + 1\right) + 1 \right]\nonumber
\end{eqnarray}
and using the properties satisfied by $G$ we conclude that
\begin{eqnarray}\label{fiduaslinha}
\varphi''(t) \leq 0.
\end{eqnarray}
It is clear that $g_{k}$ attains its supremum at some point $x_{k} \in M$. This gives the desired sequence $x_k$. It follows directly from definition of $g_{k} $ that
\begin{eqnarray*}
\grad g_{k}(x) = \grad u(x) - \varepsilon_{k}\varphi'(\gamma(x))\grad \gamma(x).
\end{eqnarray*}
In particular, at the points $x_{k}$ we obtain
\begin{eqnarray}\label{ig1}
\vert \grad u\vert(x_k) = \varepsilon_{k}\varphi'(\gamma(x_{k}))\vert \grad \gamma \vert(x_{k}).
\end{eqnarray}
Using $h2)$  in the above equality we have
\begin{eqnarray*}
\vert \grad u\vert (x_k) \leq \varepsilon_{k}.
\end{eqnarray*}

Computing $\Hess g_{k}(x)(v,v)$ we have
\begin{eqnarray}
\Hess g_{k}(x)(v,v) &=& \Hess u(x) (v,v) - \varepsilon_{k}\varphi'(\gamma(x))\Hess \gamma(x) (v,v) \nonumber \\
 && \nonumber \\ && -  \varepsilon_{k}\varphi''(\gamma(x))\langle \grad \gamma(x),v \rangle^2   \\
 &&\nonumber\\
&\geq & \Hess u(x) (v,v) - \varepsilon_{k}\varphi'(\gamma(x))\Hess \gamma(x) (v,v) \nonumber
\end{eqnarray}
for all $ v\in T_{x}M $. Using the fact that $x_{k}$ is a maximum point of $g_{k}$, the hypothesis $h4)$ and the expression for $\varphi'$, we get
\begin{eqnarray}
\Hess u(x_k) (v,v) \leq \varepsilon_{k}\varphi'(\gamma(x_k))\Hess \gamma(x_k) (v,v) \leq \varepsilon_{k}\langle v,v\rangle .
\end{eqnarray}

Finally, if assume that $h3)$ holds, we obtain
\begin{eqnarray*}
\Delta u(x_k) \leq \varepsilon_{k}\varphi'(\gamma(x_k))\Delta \gamma(x_k) \leq \varepsilon_{k}.
\end{eqnarray*}
To finish the proof of Theorem \ref{Principal Theorem} we need to show that if $u^{\ast}=\sup_{M}u<\infty$ then $u(x_{k})\to u^{\ast}$. To do that, we follow Pigola-Rigoli-Setti  closely in \cite{prs-memoirs} and observe that for any fixed $ j \in \mathbb{N} $, there is a  $ y \in M $ such that
\begin{eqnarray}\label{sup}
u(y) > \sup u - 1/2j .
\end{eqnarray}
Since $ g_k $ has a maximum at $ x_k $ we have
\begin{eqnarray*}
u(x_k) - \varepsilon_{k}\log \left(\int_{0}^{\gamma(x_k)}\frac{ds}{\sqrt{G(s)}}+1\right) = g_{k}(x_k) \geq g_{k}(y) = u(y) - \varepsilon_{k}\log \left(\int_{0}^{\gamma(y)}\frac{ds}{\sqrt{G(s)}}+1\right).
\end{eqnarray*}
Therefore, (using (\ref{sup}))
\begin{eqnarray}\label{int0}
u(x_k) > \sup u - \frac{1}{2j} - \varepsilon_{k}\log \left(\int_{0}^{\gamma(y)}\frac{ds}{\sqrt{G(s)}}+1\right).
\end{eqnarray}
Choosing  $k=k(j)=k_{j}$ sufficiently large  such that
\begin{eqnarray}\label{sup2}
\varepsilon_{k_{j}}\log \left(\int_{0}^{\gamma(y)}\frac{ds}{\sqrt{G(s)}}+1\right) < \frac{1}{2j},
\end{eqnarray}
it follows from  ({\ref{int0}}) and ({\ref{sup2}})  that
\begin{eqnarray}\label{int}
u(x_{k_j}) > \sup u - \frac{1}{2j} - \frac{1}{2j} = \sup u - \frac{1}{j}.
\end{eqnarray}
Therefore $\displaystyle\lim_{j \rightarrow +\infty} u(x_{k_j}) = \sup_{M} u $ and this finishes the proof of Theorem \ref{Principal Theorem}.

\begin{remark}
Let $\mathcal{G}$ and $\overline{\mathcal{G}}$ be the classes of Riemannian manifolds satisfying respectively the hypotheses of the Theorem \ref{r3} and Theorem \ref{Principal Theorem}. Then we have that $\mathcal{G} \subset \overline{\mathcal{G}}$.  Hence, the Theorem \ref{Principal Theorem} implies the Pigola-Rigoli-Setti's Theorem (Thm. \ref{r3}).
\begin{proof}Given $M \in \mathcal{G}$, observe that the hypothesis $$\displaystyle{\limsup_{t \rightarrow +\infty}\left[\frac{tG(\sqrt{t})}{G(t)}\right] = D < +\infty} $$ implies the existence of  $s_{0} \in \mathbb{R}$ such that
\begin{eqnarray}
\sup\left\{\frac{tG(\sqrt{t})}{G(t)}, t \geq s_{0} \right\} < D + 1.
\end{eqnarray}
Thus for all $t \geq s_{0}$ we have that
\begin{eqnarray}
\frac{tG(\sqrt{t})}{G(t)} < D + 1 & \Leftrightarrow & tG(\sqrt{t})< (D+1)G(t)
\end{eqnarray}
whence $ \sqrt{tG(\sqrt{t})} < \sqrt{(D+1)G(t)}$. In particular
\begin{eqnarray*}
A\sqrt{\gamma G(\sqrt{\gamma})} < A\sqrt{(D+1)G(\gamma)} \leq A\sqrt{(D+1)}\sqrt{G(\gamma)}\left(\int_{0}^{\gamma}\frac{ds}{\sqrt{G(s)}}+1\right).
\end{eqnarray*}
Finally, we refer to \cite[p.10]{prs-memoirs}  for a proof that any $M \in \mathcal{G}$ is a complete manifold.
\end{proof}
\end{remark}

\begin{remark}\label{remark3} 
Estimates placed the items $ h2)$, $h3)$ e $h4)$ can be exchanged for 
\begin{eqnarray}\label{condfi}
\prod_{j=1}^{k}\left[\ln^{(j)}\left(\int_{0}^{t}\frac{ds}{\sqrt{G(s)}} + 1\right)+1\right]\left(\int_{0}^{t}\frac{ds}{\sqrt{G(s)}} + 1\right)\sqrt{G(t)},
\end{eqnarray}
where $ \ln^{(j)} $ is the $j$-th iterated logarithm and $ k \in \mathbb{N} $.
\begin{proof} Indeed, note that this estimate is simply the inverse of the first derivative of the auxiliary function $ \varphi $, used in the statement of Main Theorem \ref{Principal Theorem}. Thus, we redefine the function $ \varphi $ by
\begin{eqnarray*}
\varphi(t) = \ln^{(k+1)}\left(\int_{0}^{t}\frac{ds}{\sqrt{G(s)}} + 1\right).
\end{eqnarray*}
Hence, we obtain by deriving
\begin{eqnarray*}
\varphi^{'}(t)= \left\{\prod_{j=1}^{k}\left[\ln^{(j)}\left(\int_{0}^{t}\frac{ds}{\sqrt{G(s)}} + 1\right)+1\right]\left(\int_{0}^{t}\frac{ds}{\sqrt{G(s)}} + 1\right)\sqrt{G(t)}\right\}^{-1}
\end{eqnarray*}
and
\begin{eqnarray*}
\varphi^{''}(t) &=& -\left\{\prod_{j=1}^{k}\left[\ln^{(j)}\left(\int_{0}^{t}\frac{ds}{\sqrt{G(s)}} + 1\right)+1\right]\left(\int_{0}^{t}\frac{ds}{\sqrt{G(s)}} + 1\right)\sqrt{G(t)}\right\}^{-2}\times \nonumber \\
& &\times\left\{\prod_{i=1}^{k-1}\prod_{j=1}^{i}\left[\ln^{(j)}\left(\int_{0}^{t}\frac{ds}{\sqrt{G(s)}} + 1\right)+1\right]^{-1}\left[\left(\int_{0}^{t}\frac{ds}{\sqrt{G(s)}} + 1\right)\sqrt{G(t)}\right]^{-k-1} + \right. \nonumber \\
& &\left. + \prod_{j=1}^{k}\left[\ln^{(j)}\left(\int_{0}^{t}\frac{ds}{\sqrt{G(s)}} + 1\right)+1\right]\left[\frac{G'(t)}{2\sqrt{G(t)}}\left(\int_{0}^{t}\frac{ds}{\sqrt{G(s)}} + 1\right) + 1 \right]\right\} \nonumber \\
&\leq & 0.
\end{eqnarray*}

Therefore the function $ \varphi $ satisfies the conditions necessary to prove the Main Theorem \ref{Principal Theorem}.
\end{proof}
\end{remark}

\begin{corollary}\label{corollarysectionalcondition}
Let $ M $ be a complete, noncompact, Riemannian manifold with Ricci curvature satisfying
\begin{eqnarray*}
Ric(x) \geq - B^{2}\prod_{j=1}^{k}\left[\ln^{(j)}\left(\int_{0}^{\rho(x)}\frac{ds}{\sqrt{G(s)}} + 1\right)+1\right]^{2}\left(\int_{0}^{\rho(x)}\frac{ds}{\sqrt{G(s)}} + 1\right)^{2}G(\rho(x)),
\end{eqnarray*}
for $ \rho (x)\gg 1 $ where $G\in C^{\infty}([0,+\infty)) $ satisfies \eqref{eq6}, $\rho (x)={\rm dist}_{M}(x_{0}, x)$, $B\in \mathbb{R}$. Then $\gamma=\rho$ satisfies $h1-h3$. Therefore the Omori-Yau maximum principle holds on $M$ for the Laplacian by Main Theorem \ref{Principal Theorem}. Similarly, if we assume that the radial sectional curvature satisfies the above inequality, then the Omori-Yau maximum principle holds on $M$ for the Hessian.
\end{corollary}

\vspace{.2mm}

\section{Weighted Riemannian manifolds} A weighted  manifold $(M,g,\mu_{f})$, shortly denoted by $(M,\mu_{f})$, is a Riemannian manifold  $(M, g)$  endowed with a measure $\mu_{f}=e^{-f}\nu$, where  $f\colon M\to \mathbb{R}$ is a smooth function and $\nu=\sqrt{\det g}\,dx^{1}\ldots x^{n}$ is the Riemannian density.
The associated Laplace-Betrami operator $\triangle_{f}$  is defined by $$\triangle_{f} :=e^{f}\diver (e^{-f}\grad ).$$ It is natural to extend the results above to the weighted Laplacian.
 A. Borbely, \cite{borbely} \cite{borbely1} proved a nice extension of Pigola-Rigoli-Setti's version \cite{prs-memoirs} of the Omori-Yau maximum principle for the Laplacian. Borbely's version has been  extended to the weighted Laplacian or even to more general operators  by many authors. Some in the weak form of the maximum principle others in the strong form of the maximum principle.  For instance,  Bessa, Pigola and Setti in
 \cite[Thm 9]{bps-revista}, Pigola Rigoli and Setti \cite{prs-jfa},  Mari, Rigoli and Setti  in \cite{mari-rigoli-setti}, by  Pigola, Rigoli, Rimoldi and Setti in \cite{pigola-rigoli-rimoldi} and by Mastrolia, Rigoli and Rimoldi  in \cite{matrolia-rigoli-rimoldi}.  For the Laplace operator we can resume what they proved as
\begin{theorem}[Borb\'{e}ly, Bessa, Mari, Mastrolia, Pigola, Rigoli, Rimoldi, Setti]Let $(M,\mu_{f})$ be a complete weighted manifold and assume that there exists a non-negative $C^{2}$-function $\gamma$ satisfying the following conditions.
\begin{itemize}\item[a.] $\gamma (x)\to + \infty$ as $x\to \infty$.
\item[]
\item[b.]$\exists A>0$ such that  $\vert \grad \gamma \vert <A$ off a compact set.
\item[]
\item[c.]$\exists B>$ such that $\triangle_{f} \gamma \leq B\cdot G(\gamma )$ off a compact set.
\end{itemize}Where $G\in C^{\infty}([0, \infty))$ satisfying 
\begin{equation}\label{eqthm5}
\begin{array}{lllll} G(0)>0,&& G'(t)\geq 0 \,\, in\,\, [0, \infty),&& G(t)^{-1}\not \in L^{1}([0, \infty).\end{array}\end{equation}Then the Omori-Yau maximum principle for $\triangle_{f}$ holds on $M$.\label{bps-revista}
\end{theorem}

The main result of this section is the following  extension of Theorem \ref{bps-revista}.

\begin{theorem}Let $(M,\mu_{f})$ be a complete weighted manifold and assume that there exists a non-negative $C^{2}(M)$-function $\gamma$ satisfying the following conditions:
\begin{itemize}\item[a.] $\gamma (x)\to + \infty$ as $x\to \infty$;
\item[]
\item[b.]$\exists A>0$ such that  $\vert \grad \gamma \vert <A \cdot \sqrt{G(\gamma)}\left(\int_{0}^{\gamma}\frac{ds}{\sqrt{G(\gamma)}}+1\right) $ off a compact set;
\item[]
\item[c.]$\exists B>$ such that $\triangle_{f} \gamma \leq B \cdot \sqrt{G(\gamma)}\left(\int_{0}^{\gamma}\frac{ds}{\sqrt{G(\gamma)}}+1\right) $ off a compact set; \\
where $G\in C^{\infty}([0, \infty))$ satisfies 
\begin{equation}\label{eqthm6}
\begin{array}{lllll} G(0)>0,&& G'(t)\geq 0 \,\, in\,\, [0, \infty),&& G(t)^{-1/2}\not \in L^{1}([0, \infty).\end{array}\end{equation}
\end{itemize}Then the Omori-Yau maximum principle for $\triangle_{f}$ holds on $M$.
\end{theorem}

\begin{remark}Replacing the bound $\triangle_{f} \gamma <A\cdot  G(r)$ in Theorem \ref{bps-revista}, with  $G$ satisfying \eqref{eqthm5},   by $\triangle_{f} \gamma \leq B \cdot \sqrt{\tilde{G}(\gamma)}\left(\int_{0}^{\gamma}\frac{ds}{\sqrt{\tilde{G}(\gamma)}}+1\right) $, with $\tilde{G}$ satisfying  \eqref{eqthm6} does not amount to a weaker condition. In fact, if $\tilde{G}$ satisfies \eqref{eqthm6} then $G(r)= \sqrt{\tilde{G}}$ satisfies \eqref{eqthm5}.
\end{remark}

\begin{remark}We remark that   if $u\in C^{2}(M)$ with $
\displaystyle{\lim_{x \to \infty}\frac{u(x)}{\varphi(\gamma(x))}} = 0 ,
$ then there exist a sequence ${x_{k}} \in M$, $k \in \mathbb{N} $ such that
$ \displaystyle{\triangle_{f} u(x_{k}) < \frac{1}{k}} $.
\end{remark}

\begin{proof}
The proof follows closely the  proof of the Main Theorem \ref{Principal Theorem}. We need only to adapt the part of the proof that treats with the $\triangle$ to $\triangle_{f}$. As in there we defined the functions sequence $ g_k $ and observe that
\begin{eqnarray*}
\triangle_{f}\varphi(\gamma)(x) &=&  \varphi'(\gamma(x))\triangle_{f}\gamma(x) + \varphi''(\gamma(x))\vert \grad \gamma\vert^2(x) .
\end{eqnarray*}
Note that in the maximum points $ x_k $ of the functions $ g_k $, we have $ \triangle_{f}g_{k}(x_{k}) \leq 0 $. Thus,
\begin{eqnarray*}
0 &\geq & \triangle_{f}g_{k}(x_k) = \triangle_{f}u(x_k) - \varepsilon_{k}\triangle_{f}\varphi(\gamma(x_{k})) \\
&=& \triangle_{f}u(x_k) - \varepsilon_{k}\left[\varphi'(\gamma(x_{k})\triangle_{f}\gamma(x_k) + \varphi''(\gamma(x_k))\vert\grad \gamma(x_k)\vert^2\right],
\end{eqnarray*}
which implies in
\begin{eqnarray*}
\triangle_{f}u(x_k) &\leq & \varepsilon_{k}[\varphi'(\gamma(x_k))\triangle_{f}\gamma(x_k) + \varphi''(\gamma(x_k))\vert\grad \gamma(x_k)\vert^2] \\
&\leq & \varepsilon_{k}\varphi'(\gamma(x_k))\triangle_{f}\gamma(x_k).
\end{eqnarray*}
In the last inequality we used that $\varphi''\leq 0$. Then the we proceed as in Theorem \ref{Principal Theorem} to finish the proof.
\end{proof}

\vspace{.2mm}

\section{Geometric Applications}In a beautiful paper \cite{jorge-koutrofiotis}, Jorge and Koutrofiotis  applied   Omori's Theorem \eqref{r1} to give curvature estimates for bounded submanifolds with scalar curvature bounded below, extending various non-immersability results. Their result was extended by Pigola, Rigoli and Setti in \cite{prs-memoirs} as an application of their generalized version of the Omori-Yau maximum principle.

Recently, L. Alias, G. P. Bessa and J. F. Montenegro in \cite{alias-bessa-montenegro} proved a  version of Jorge-Koutrofiotis Theorem for cylindrically bounded submanifolds, recalling that an isometric immersion $\varphi \colon M\hookrightarrow N\times \mathbb{R}^{\ell}$ is said to be cylindrically bounded if $\varphi (M)\subset B_{N}(r)\times \mathbb{R}^{\ell}$, where $B_{N}(r)$ is a geodesic ball in $N$ of radius $r>0$.

\begin{theorem}
\label{1thmMain}$[$Alias-Bessa-Montenegro$]$
Let $M$ and $N$  be complete Riemannian manifolds of dimension $m$ and $n-\ell$ respectively, satisfying
$n+\ell\leq 2m-1$. Let $\varphi:M^m\rightarrow N^{n-\ell}\times\R^{\ell}$ be a isometric immersion with  $\varphi(M)\subset B_N(r)\times\R^{\ell}$.
Assume that the radial sectional curvature $K_N^{\mathrm{rad}}$ along the radial geodesics issuing from $p$
satisfies $K_N^{\mathrm{rad}}\leq b$ in $B_N(r)$ and $0<r<\min\{\mathrm{inj}_N(p),\pi/2\sqrt{b}\}$, where we replace
$\pi/2\sqrt{b}$ by $+\infty$ if $b\leq 0$. Suppose that the immersion $\varphi $ is proper and
\begin{equation}
\sup_{\varphi^{-1}(B_{N}(r)\times \partial B_{\mathbb{R}^{\ell}}(t))}\Vert\alpha\Vert\leq \sigma(t),
\label{growth}
\end{equation}
where $\alpha$ is the second fundamental form of the immersion and $\sigma:[0,+\infty)\rightarrow\mathbb{R}$ is a
positive  function satisfying $\int_0^{+\infty}1/\sigma=+\infty$, then the sectional curvature of $M$ has the following lower bound
\begin{equation}
\sup_{M}K_{M}\geq C_{b}^{2}(r)+\inf_{B_{N}(r)}K_{N},
\label{eq-ABM}
\end{equation} with
\[
C_b(t)=
\begin{cases}
\sqrt{b}\cot(\sqrt{b}\,t) & \mbox{\rm if $b>0$ and $0<t<\pi/2\sqrt{b}$}\\
1/t & \mbox{\rm if $b=0$ and $t>0$}\\
\sqrt{-b}\coth(\sqrt{-b}\,t) & \mbox{\rm if $b<0$ and $t>0$}.
\end{cases}
\]
\end{theorem}In this section our main result is the following generalization of Theorem \ref{1thmMain}. We prove the following result.
\begin{theorem}\label{thmMain}
Let $M$ and $N$  be complete Riemannian manifolds of dimension $m$ and $n-\ell$ respectively, satisfying
$n+\ell\leq 2m-1$. Let $\varphi:M^m\rightarrow N^{n-\ell}\times\R^{\ell}$ be a proper isometric immersion with  $\varphi(M)\subset B_N(r)\times\R^{\ell}$.
Assume that the radial sectional curvature $K_N^{\mathrm{rad}}$ along the radial geodesics issuing from $p$
satisfies $K_N^{\mathrm{rad}}\leq b$ in $B_N(r)$ and $0<r<\min\{\mathrm{inj}_N(p),\pi/2\sqrt{b}\}$, where we replace
$\pi/2\sqrt{b}$ by $+\infty$ if $b\leq 0$.
 Then the sectional curvature of $M$ has the following lower bound
\begin{equation}
\sup_{M}K_{M}\geq C_{b}^{2}(r)+\inf_{B_{N}(r)}K_{N}.
\label{eq-ABM}
\end{equation}
\end{theorem}
\begin{proof}Let $ g : N\times\mathbb{R}^{\ell} \rightarrow \mathbb{R} $ be given $ g(z,y) = \phi_{b}(\rho_{N}(z)) $, where $\phi_{b}$ is given by
\begin{eqnarray}
\phi_{b}(t) =
\left\{\begin{array}{ll}
1 - \cos(\sqrt{b}\,t) & \mbox{if} \ b > 0 \ \mbox{and} \ 0 < t < \pi/2\sqrt{b} \\
t^2 & \mbox{if} \ b = 0 \ \mbox{and} \ t > 0 \\
\cosh(\sqrt{-b}\,t) & \mbox{if} \ b < 0 \ \mbox{and} \ t > 0
\end{array}\right.
\end{eqnarray}
and $ \rho_{N}(z) = dist_{N}(p,z) $.  Consider $ f : M \rightarrow \mathbb{R} $, $ f = g\circ\varphi $ and let $ \pi_{N} : N\times\mathbb{R}^{\ell} \rightarrow N $ be the projection on the factor $ N $. Since $ \pi_{N}(\varphi(M)) \subset B_{N}(r) $, we have that $ f^{*} = \sup_{M}f \leq \phi_{b}(r) < + \infty $. Define for each $ k \in \mathbb{N} $, the function $ g_{k} : M \rightarrow \mathbb{R} $ be given
\begin{eqnarray}
g_{k}(x) = f(x) - \varepsilon_{k}\psi(\rho_{\mathbb{R}^{\ell}}(y(x))),
\end{eqnarray}
where $ \psi : \mathbb{R} \rightarrow [0,+\infty) $ is given by $\psi(t) = \log(\log(t+1)+1) $, $\rho_{\mathbb{R}^{\ell}}(y)={\rm dist}_{\mathbb{R}^{\ell}}(0, y)$, $\varepsilon_k\to 0^+$ as $k\to \infty$
and $ y(x) = \pi_{\mathbb{R}^{\ell}}(\varphi(x)) $.

Since the immersion $ \varphi $ is proper, if $ x \rightarrow \infty $ in $M$  then $\varphi (x)\to \infty$ in $B_{N}(r)\times\mathbb{R}^{\ell}$, thus $y(x)\to \infty$ in $\mathbb{R}^{\ell}$ and  $ \psi(\rho_{\mathbb{R}^{\ell}}(y(x))) \rightarrow +\infty $. Therefore $g_{k}$ reach its maximum at a point $x_{k}\in M$. This forms a sequence $\{ x_k \}\subset M $ such that $ g_{k}(x_k) = \sup_{M}g_k $.  There are two cases to consider:
 \begin{itemize}\item[1.]$x_k\to \infty$ in $M$ as $k\to +\infty$. \item[2.] $x_k$ stays in a bounded subset of $M$. \end{itemize}
  Let us consider the case 1. i.e. $x_k\to \infty$ in $M$ as $k\to +\infty$.
Since $x_{k}$ is a point of maximum for $g_{k}$ we have that $\Hess_{M}g_{k}(X,X)\leq 0$ for all $X\in T_{x_{k}}M$. This implies that \begin{equation}\Hess_{M}f(x_{k})(X,X)\leq \varepsilon_{k}\Hess_{M}\psi\circ \rho_{\mathbb{R}^{\ell}}\circ y(x_{k})(X,X),\,\,\,  X\in T_{x_{k}}M.\label{eqHessiano}\end{equation} 

 \noindent First we will compute the right hand side of \eqref{eqHessiano}.  We have then
\begin{eqnarray}\Hess_{M}\psi\circ \rho_{\mathbb{R}^{\ell}}\circ y(x_{k})(X,X)&=&\Hess_{N\times \mathbb{R}^{\ell}}\psi\circ \rho_{\mathbb{R}^{\ell}}\circ y(x_{k})(X,X)\nonumber \\
&&\label{eqHessiano2}\\
&&+ \langle \grad_{N\times \mathbb{R}^{\ell}} \psi\circ \rho_{\mathbb{R}^{\ell}}\circ y(x_{k}), \alpha_{M}(X,X) \rangle ,  \nonumber
\end{eqnarray}
where $\alpha$ is the second fundamental form of the immersion $\varphi$, see \cite{jorge-koutrofiotis}.
\vspace{2mm}

\noindent
Setting  $y_{k}=\pi_{\mathbb{R}^{\ell}}(\varphi (x_{k}))$ and $t_{k}=\rho_{\mathbb{R}^{\ell}}(y_k)$ we have 
\begin{eqnarray}
\Hess_{N\times \mathbb{R}^{\ell}}\psi\circ \rho_{\mathbb{R}^{\ell}}\circ y(x_{k})(X,X) &=&  \psi''(t_k)\vert X^{\mathbb{R}^{\ell}}\vert^{2}+\psi'(t_k)\Hess_{\mathbb{R}^{\ell}} \rho_{\mathbb{R}^{\ell}}(y_k)(X,X)\nonumber\\
&&\nonumber\label{eqHessiano3} \\
&=&   \psi''(t_k)\vert X^{\mathbb{R}^{\ell}}\vert^{2}+ \frac{\vert X^{N}\vert^{2}}{t_k(t_k+1)(\log (t_k+1)+1)}\\
&&\nonumber \\ 
&\leq & \frac{\vert X\vert^{2}}{t_k(t_k+1)(\log (t_k+1)+1)}\nonumber 
\end{eqnarray}Since $\psi'' \leq 0$.
 Here $X^{\mathbb{R}^{\ell}}=d\pi_{\mathbb{R}^{\ell}}X$ and  $X^{N}= d\pi_{N}X$, where $\pi_{\mathbb{R}^{\ell}}\colon N\times \mathbb{R}^{\ell}\to \mathbb{R}^{\ell} $, $\pi_{N}\colon N\times \mathbb{R}^{\ell}\to N $ are standard projections.
 
 We also have 
 \begin{eqnarray}\langle \grad_{N\times \mathbb{R}^{\ell}} \psi\circ \rho_{\mathbb{R}^{\ell}}\circ y(x_{k}), \alpha_{M}(X,X)\rangle &=& \psi'(t_k)\langle \grad \rho_{\mathbb{R}^{\ell}}(y_k), \alpha (X, X)\rangle  \nonumber\\ 
 && \nonumber \\
 &\leq & \frac{1}{(t_k+1)(\log (t_k+1)+1)}\vert \alpha(X,X)\vert \label{eqHessiano4} \end{eqnarray}

 From \eqref{eqHessiano3} and \eqref{eqHessiano4} we have that 
 \begin{eqnarray}\Hess_{M}\psi\circ \rho_{\mathbb{R}^{\ell}}\circ y(x_{k})(X,X)&\leq & \frac{1+ \vert \alpha(X,X)\vert  }{(t_k+1)(\log (t_k+1)+1)}\vert X\vert^{2}\label{eqHessiano5}
 \end{eqnarray}And from \eqref{eqHessiano} and \eqref{eqHessiano5} we get \begin{eqnarray}\Hess_{M}f(x_{k})(X,X)& \leq & \frac{\varepsilon_{k}(1+ \vert \alpha(X,X)\vert  )}{(t_k+1)(\log (t_k+1)+1)}\vert X\vert^{2}\label{eqHessianoI}
 \end{eqnarray}

Now,  we will compute the left hand side of \eqref{eqHessiano}.
\begin{eqnarray}\Hess_{M}f(x_{k})(X,X) &=& \Hess_{N\times \mathbb{R}^{\ell}} g(\varphi (x))(X,X) + \langle \grad g, \alpha (X,X)\rangle\label{eqHessiano6}
\end{eqnarray}Recalling that $f=g\circ\varphi$ and $g$ is given by  $ g(z,y) = \phi_{b}(\rho_{N}(z)) $, where $\phi_{b}$ is given by
\begin{eqnarray}
\phi_{b}(t) =
\left\{\begin{array}{ll}
1 - \cos(\sqrt{b}\,t) & \mbox{if} \ b > 0 \ \mbox{and} \ 0 < t < \pi/2\sqrt{b} \\
t^2 & \mbox{if} \ b = 0 \ \mbox{and} \ t > 0 \\
\cosh(\sqrt{-b}\,t) & \mbox{if} \ b < 0 \ \mbox{and} \ t > 0.
\end{array}\right.
\end{eqnarray}
and $ \rho_{N}(z) = dist_{N}(p,z) $.  Let us consider an orthonormal basis $$\{\stackrel{\in TN}{\overbrace{ \grad \rho_{N}, \partial/\partial \theta_{1}, \ldots, \partial/\partial \theta_{n-\ell-1}}}, \stackrel{\in T\mathbb{R}^{\ell}}{\overbrace{ \partial/\partial \gamma_1, \ldots, \partial/\partial \gamma_\ell }}\}$$ for $T_{\varphi (x_k)}(N\times \mathbb{R}^{\ell})$. Thus if $X\in T_{x_k}M$, $\vert X\vert =1$, we can decompose $$X=a\cdot \grad \rho_{N}+\sum_{j=1}^{n-\ell-1}b_{j}\cdot \partial/\partial \theta_{j} + \sum_{i=1}^{\ell}c_{i}\cdot \partial/\partial \gamma_i$$  with $a^2+\sum_{j=1}^{n-\ell-1} b_{j}^{2}+ \sum_{i=1}^{\ell} c_{i}^{2}=1$. Letting $z_k=\pi_{N}(\varphi (x_k))$ and $s_k=\rho_{N}(z_k)$. Having set that we have that the first term of the right hand side of \eqref{eqHessiano6}
\begin{eqnarray} \Hess_{N\times \mathbb{R}^{\ell}} g(\varphi (x))(X,X)&=&\!\phi_{b}''(s_k)\cdot   a^2+ \phi_{b}'(s_k)\sum_{j=1}^{n-\ell-1}\!b_j^{2}\cdot \Hess \rho_{N}(z_k)(\frac{\partial}{\partial \theta_{j}},\frac{\partial}{\partial \theta_{j}})\nonumber\\
&&\nonumber \\
&\geq &\!\phi_{b}''(s_k)\cdot   a^2+ \phi_{b}'(s_k)\sum_{j=1}^{n-\ell-1}\!b_j^{2}\cdot C_{b}(s_k)\nonumber\\
&&\nonumber \\
&=& \phi_{b}''(s_k)\cdot a^2+   (1-a^2 - \sum_{i=1}^{\ell}\!c_i^{2})\cdot \phi_{b}'(s_k)\cdot C_{b}(s_k)\nonumber \\
&& \nonumber \\
&=& \left[\stackrel{\equiv 0}{(\overbrace{\phi_{b}''-C_{b}\cdot \phi_{b}'})}a^2 + (1- \sum_{i=1}^{\ell}\!c_i^{2})\cdot \phi_{b}'\cdot C_{b}\right](s_k) \nonumber\\
&&\nonumber \\
&=& (1- \sum_{i=1}^{\ell}\!c_i^{2})\cdot \phi_{b}'(s_k)\cdot C_{b}(s_k).\nonumber
\end{eqnarray}

Thus \begin{equation}\label{eqHessiano7}\Hess_{N\times \mathbb{R}^{\ell}} g(\varphi (x))(X,X)\geq (1- \sum_{i=1}^{\ell}\!c_i^{2})\cdot \phi_{b}'(s_k)\cdot C_{b}(s_k)\cdot \vert X\vert^{2}.\end{equation}

We used  above two facts.  The first was $\Hess \rho_{N}(z_k)(\frac{\partial}{\partial \theta_{j}},\frac{\partial}{\partial \theta_{j}})\geq C_{b}(s_{k})$ yielded by  the Hessian Comparison Theorem, Thm. \ref{thmHessiano}. Recall that the radial sectional curvature of $N$ along the geodesics issuing from the center of  the ball $B_{N}(r)$ is bounded above $K_{N}^{rad}\leq b$, see the hypotheses of Theorem \ref{thmMain}. We state  the Hessian Comparison Theorem for sake of completeness.
 The second fact is the $\phi_{b}$ satisfies the following equation $ \phi_{b}''(t)-C_{b}(t)\phi_{b}'(t)=0$, $C_{b}$ given below in \eqref{eqCb} .
\begin{theorem}[Hessian Comparison Theorem] \ Let $ M $ be a Riemannian manifold and $ x_{0},x_{1} \in M $ be such that there is a minimizing unit speed geodesic $ \gamma $ joining $x_0 $ and $ x_1 $ and let $\rho(x) = dist(x_{0},x) $ be the distance function to $ x_0 $. Let $ K_{\gamma} \leq b $ be the radial sectional curvatures of $ M $ along $ \gamma $. If $ b > 0 $ assume $ \rho(x_{1}) < \pi/2\sqrt{b} $. Then, we have $ \Hess\rho(x)(\dot{\gamma},\dot{\gamma}) = 0 $ and
\begin{eqnarray}
\Hess \rho(x)(X,X) \geq C_{b}(\rho(x))\vert X\vert^2
\end{eqnarray}
where $ X \in T_{x}M $ is perpendicular to $ \dot{\gamma}(\rho(x)) $ and
\begin{eqnarray}\label{eqCb}
C_{b}(t) =
\left\{\begin{array}{ll}
\sqrt{b}\cot(\sqrt{b}t) & \mbox{if} \ b > 0 \ \mbox{and} \ 0 < t < \pi/2\sqrt{b} \\
1/t & \mbox{if} \ b = 0 \ \mbox{and} \ t > 0 \\
\sqrt{-b}\coth(\sqrt{-b}t) & \mbox{if} \ b < 0 \ \mbox{and} \ t > 0.
\end{array}\right.
\end{eqnarray}\label{thmHessiano}
\end{theorem}

The second term of the right hand side of \eqref{eqHessiano6} is the following, if $\vert X\vert =1$.
\begin{eqnarray} \langle \grad g, \alpha (X,X)\rangle & = & \phi_{b}'(s_k)\langle \grad \rho_{N}(y_k), \alpha (X,X)\rangle\nonumber\\ 
&& \label{eqHessiano8}\\
&\geq & -\phi_{b}'(s_k)\vert \alpha(X,X) \vert\nonumber
\end{eqnarray}

Therefore from \eqref{eqHessiano6}, \eqref{eqHessiano7}, \eqref{eqHessiano8} we have that

\begin{equation}\Hess_{M}f(x_{k})(X,X) \geq \left[(1- \sum_{i=1}^{\ell}\!c_i^{2})\cdot  C_{b}(s_k) -  \vert \alpha (X/\vert X\vert,X/\vert X\vert) \vert\right]\phi_{b}'(s_k)\vert X\vert^{2}\label{eqHessiano9}\end{equation}

Recall that we have an isometric immersion $\varphi \colon M^{m}\hookrightarrow N^{n-\ell}\times \mathbb{R}^{\ell}$, where $n+\ell\leq 2m-1$. This dimensional restriction implies that $m\geq \ell + 2$. Therefore, for every $x\in M$ there exists a  sub-space $V_x\subset T_xM$ with ${\rm dim}(V_x)\geq (m-\ell)\geq  2$ such that $V\perp T\mathbb{R}^{\ell}$. If we take any $X\in V_{x_k}\subset T_{x_k}M$, $\vert X\vert=1$ we have by \eqref{eqHessianoI}, \eqref{eqHessiano9} that

\begin{equation}\frac{\varepsilon_{k}(1+ \vert \alpha(X,X)\vert  )}{(t_k+1)(\log (t_k+1)+1)}\geq \Hess_{M}f(x_{k})(X,X) \geq \left[  C_{b}(s_k) -  \vert \alpha (X,X) \vert\right]\phi_{b}'(s_k)\label{eqHessiano10}\end{equation} Recall that $t_k = \rho_{\mathbb{R}^{\ell}}(\pi_{\mathbb{R}^{\ell}}(\varphi (x_{k})))\to \infty$, $s_{k}=\rho_{N}(\pi_{N}(\varphi (x_k)))<r$ and $\varepsilon_k\to 0^+$ as $k\to \infty$.

From \eqref{eqHessiano10} we have that 
\begin{eqnarray*}\vert \alpha(x_k) (X,X)\vert \left[ \frac{\varepsilon_k}{(t_k +1)(\log (t_k+1)+1)} + \phi_{b}'(s_k)\right]&\geq&  C_{b}(s_k)\phi_{b}'(s_k)\\
&& \\
&-& \frac{\varepsilon_k}{(t_k+1)(\log (t_k+1)+1)} .
\end{eqnarray*}
Thus  for  large $k$ and for all $0\neq X\in  V_{x_k}\subset T_{x_k}M$ we have that 
$$\vert\alpha (x_k)(X,X)\vert \geq  \left[C_b(s_k)+ \stackrel{\to 0 }{\overbrace{\delta (k, \varepsilon_k, t_k)}}\right]\vert X\vert^{2}.$$ We will need the following lemma known as Otsuki's Lemma. \begin{lemma}[Otsuki] Let $ \beta \colon \mathbb{R}^{k}\times\mathbb{R}^k \rightarrow \mathbb{R}^d $, $ d \leq k-1 $, be a symmetric bilinear form satisfying $ \beta(X,X)\neq 0 $ for $ X\neq 0 $. Then there exists linearly independent vectors $ X,Y $ such that $\beta(X,X)=\beta(Y,Y) $ and $ \beta(X,Y)=0 $.
\end{lemma} Observe that we just showed that for all $0\neq X\in  V_{x_k}\subset T_{x_k}M$ we have that $\vert\alpha (x_k)(X,X)\vert > 0. $ Moreover, ${\rm dim}(V_{x_k})=m-\ell \leq n-m$ by hypothesis. Applying Otsuki's Lemma to $\alpha (x_k)\colon  V_{x_k}\times V_{x_k}\to T_{x_k}M^{\perp}$ we obtain $X, Y\in V_{x_k}$, $\vert X\vert \geq \vert Y\vert \geq 1$ such that $\alpha (x_k)(X,X)= \alpha (x_k)(Y,Y)$ and $\alpha (x_k)(X,Y)=0$.

Using Gauss equation we have that 
\begin{eqnarray}\label{curvaturek}
K_{M}(x_k)(X,Y) - K_{N}(\varphi(x_k))(X,Y) &=& \frac{\langle \alpha (x_k)(X,X),\alpha (x_k)(Y,Y)\rangle - \vert \alpha (x_k)(X,Y)\vert^{2} }{\vert X\vert^{2}\vert Y\vert^{2} - \langle X,Y\rangle^{2}} \nonumber \\
&\geq & \frac{\vert \alpha (x_k)(X,X)\vert^{2}}{\vert X\vert^{2}\vert Y\vert^{2}} \nonumber \\
&\geq & \left(\frac{\vert \alpha (x_k)(X,X)\vert }{\vert X\vert^{2}}\right)^{2}\nonumber \\
&\geq &  \left(C_{b}(s_k)+ \stackrel{\to 0 }{\overbrace{\delta (k, \varepsilon_k, t_k)}}\right)^{2}.
\end{eqnarray}Letting $k\to\infty$ we obtain that $C_b(s_k)\to C_{b}(s^{\ast})\geq C_{b}(r)$, $s^{\ast}\leq r$ and we have that $$ \sup K_{M}-\inf K_{N}\geq C_{b}^{2}(r)$$

The  case where the sequence $ \{x_k\} \subset M $ remains in a compact set we proceed as follows.  Passing to a subsequence we have that $ x_k \rightarrow x_0 \in M $ and $ f $ attains its absolute maximum at $ x_0 $. Thus $ f(x_0)(X,X) \leq 0 $ for all $ X \in T_{x_0}M $. It using the expression on the right hand side of \eqref{eqHessiano10} we obtain for every $ X \in V_{x_0} $
\begin{eqnarray*}
0 \geq \Hess f(x_0)(X,X) \geq \phi'_{b}(s_0)\left(C_{b}(s_0)\vert X\vert ^{2} - \vert \alpha_{x_0}(X,X)\vert\right) .
\end{eqnarray*}
Hence
\begin{eqnarray*}
\vert \alpha_{x_0}(X,X)\vert  \geq C_{b}(s_0)\vert X\vert^2 .
\end{eqnarray*}

Following the step made in the statement above, we conclude that
\begin{eqnarray}
\sup_{M}K_{M} - \inf_{B_{N}(r)}K_{N} \geq C_{b}^{2}(s_0) \geq C_{b}^{2}(r) .
\end{eqnarray}
This finishes the proof of Theorem \ref{thmMain}.
\end{proof}

Following the terminology introduced in \cite{alias-bessa-montenegro-piccione}, we indroduce the next definition.

\begin{definition}
The pair of functions $ (G,\gamma) $ where $ G : [0,+\infty) \rightarrow [0,+\infty) $ and \linebreak 
$ \gamma : M \rightarrow [0,+\infty) $ form an Omori-Yau pair for the Hessian, respectively Laplacian, in M if they satisfy the conditions established in the Main Theorem \ref{Principal Theorem} (or the condition $(\ref{condfi})$ of the Remark $\ref{remark3}$).
\end{definition}
Recently, Al\'ias and  Dajczer in \cite{alias-dajczer}, studying the mean curvature estimates for cylindrically bounded submanifolds, showed that if we take a proper isometric immersion $ \varphi : M^{m} \rightarrow L^{\ell}\times_{\rho}P^n $  then the existence of a Omori-Yau pair for the Hessian in $ L^\ell $ induces an Omori-Yau pair for the Laplacian on $M^{m}$ provided the mean curvature  $ \vert H \vert $ is bounded. Recalling that $ L^{\ell}\times_{\rho}P^n $ is endowed with product metric $ ds^{2} = dg_{L} + \rho^{2}dg_{P} $, where $ dg_{L} $ and $ dg_{P} $ are the Riemannian metrics of $ L $ and $ P $ respectively.

The next proposition generalizes the essential fact in the proof of Theorem 1 in \cite{alias-dajczer}.

\begin{proposition}\label{pullbackomoriyau} Let $ \varphi : M^{m} \rightarrow L^{\ell}\times_{\rho}P^n = N^{n+l} $ be an isometric immersion where $ L^\ell $ carries an Omori-Yau pair $(G,\tilde{\gamma}) $ for the Hessian, $ \rho \in C^{\infty}(L) $ is a positive function and the function $ \vert \grad \log \rho\vert $ satisfies
\begin{eqnarray}\label{hipeta}
\vert \grad \log \rho\vert \leq \ln\left(\int_{0}^{\tilde{\gamma}}\frac{ds}{\sqrt{G(s)}} + 1\right).
\end{eqnarray}
If $ \varphi $ is proper on the first entry and
\begin{eqnarray}\label{hipcurvmedia}
\vert H\vert \leq \ln\left(\int_{0}^{\tilde{\gamma}\circ\pi_{L}}\frac{ds}{\sqrt{G(s)}}+1\right),
\end{eqnarray}then $ M^m $ has an Omori-Yau pair for the Laplacian.
Here $ \pi_{L}\colon L\times P\to L $ is the projection on factor $L$.
\end{proposition}
\begin{proof} The crux of the proof is presented in \cite{alias-dajczer} and therefore will try to follow the same notation to simplify the proof.
Suppose that $ M $ is non-compact and denote $ \varphi = (x,y) $. Define $ \Gamma(x,y) = \tilde{\gamma}(x) $ and $ \gamma = \Gamma\circ \varphi = \tilde{\gamma}(x) $. We claim that $ (G,\gamma) $ is an Omori-Yau pair for the laplacian in $ M $.

Indeed, let $ q_k \in M $ a sequence such that $ q_k \rightarrow \infty $ in $ M $ as $ k \rightarrow +\infty $. Since $ \varphi $ is proper in the first entry, we have that $ x(q_k) \rightarrow \infty $ in $ L $. Hence $ \gamma(q_k) \rightarrow \infty $ as $ k \rightarrow +\infty $, because $ \tilde{\gamma} $ also is proper. 

We have from $ \Gamma(x,y) = \tilde{\gamma}(x) $ that
\begin{eqnarray}\label{gradGama}
\grad \Gamma(x,y) = \grad \tilde{\gamma}(x) .
\end{eqnarray}
Since $ \gamma = \Gamma\circ\varphi $, we obtain
\begin{eqnarray}
\grad \Gamma(\varphi(q)) = \grad \gamma(q) + (\grad \Gamma(\varphi(q)))^{\perp} ,
\end{eqnarray}
and by the hypothesis, we have
\begin{eqnarray*}
\vert \grad \gamma \vert(q) &\leq & \vert \grad \Gamma\vert (\varphi(q)) \nonumber \\
&=& \vert \grad \tilde{\gamma}\vert(x(q)) \nonumber \\
&\leq & \sqrt{G(\gamma(q))}\left(\int_{0}^{\gamma(q)}\frac{ds}{\sqrt{G(s)}} + 1\right) ,
\end{eqnarray*}
outside a compact subset of $ M $.

Since $ \nabla^{N}_{S}T = \nabla^{L}_{S}T $, for all $ T,S \in TL $, follows from (\ref{gradGama}) that
\begin{eqnarray*}
\nabla_{T}^{N}\grad \Gamma = \nabla^{L}_{T}\grad \tilde{\gamma} .
\end{eqnarray*}
Hence,
\begin{eqnarray*}
\Hess\Gamma(T,S) = \Hess\tilde{\gamma}(T,S) \quad \mbox{and} \quad \Hess\Gamma(T,X) = 0,
\end{eqnarray*}
where $ T,S \in TL $ and $ X \in TP $. Moreover, since $ \nabla^{N}_{X}T = \nabla^{N}_{T}X = T(\eta)X $, for all $ T \in TL $, $ X \in TP $ and $ \eta = \log\rho $, we have 
\begin{eqnarray*}
\nabla^{N}_{X}\grad \Gamma = \grad \tilde{\gamma}(\eta)X .
\end{eqnarray*}
Thus,
\begin{eqnarray*}
\Hess\Gamma(X,Y) &=& \langle \nabla^{N}_{X}\grad \Gamma,Y\rangle \nonumber \\
&=& \langle \grad \tilde{\gamma}(\eta)X,Y\rangle \nonumber \\
&=& \langle \langle \grad \eta,\grad \tilde{\gamma}\rangle X,Y\rangle \nonumber \\
&=& \langle \grad \tilde{\gamma},\grad \eta\rangle\langle X,Y \rangle .
\end{eqnarray*}

For a unit vector $ e \in T_{q}M $, set $ e = e^L + e^P $, where $ e^L \in T_{x(q)}L $ and $ e^P \in T_{y(q)}P $. Then we have
\begin{eqnarray*}
\Hess\Gamma(\varphi(q))(e,e) = \Hess\tilde{\gamma}(x(q))(e^L,e^L) + \langle \grad \tilde{\gamma}(x(q)),\grad \eta(x(q))\rangle\vert e^P\vert^2 .
\end{eqnarray*}
Since $ \gamma = \Gamma\circ\varphi $, we get
\begin{eqnarray}\label{expression1}
\Hess \gamma(q)(e,e) &=& \Hess\tilde{\gamma}(x(q))(e^L,e^L) + \langle \grad \tilde{\gamma}(x(q)),\grad \eta(x(q))\rangle\vert e^P\vert^2 + \nonumber \\
&&+ \langle \grad \tilde{\gamma}(x(q)),\alpha_{q}(e,e)\rangle .
\end{eqnarray}

But,
\begin{eqnarray}\label{expressionI}
\Hess\tilde{\gamma} (.,,) \leq \sqrt{G(\tilde{\gamma})}\left(\int_{0}^{\tilde{\gamma}}\frac{ds}{\sqrt{G(s)}} + 1\right)\langle .,.\rangle ,
\end{eqnarray}
outside a compact subset in $ L $ and by the hypothesis \eqref{hipeta}
\begin{eqnarray}\label{expressionII}
\langle \grad \tilde{\gamma}(x(q)),\grad \eta(x(q))\rangle\vert e^P\vert^2 &\leq & \vert \grad \tilde{\gamma}(x(q))\vert \cdot \vert \grad \eta(x(q))\vert \nonumber \\
&\leq & \sqrt{G(\gamma)}\left(\int_{0}^{\gamma}\frac{ds}{\sqrt{G(s)}} + 1\right)\ln\left(\int_{0}^{\gamma}\frac{ds}{\sqrt{G(s)}} + 1\right).
\end{eqnarray}

Considering $(\ref{expressionI})$ and $(\ref{expressionII})$ into $(\ref{expression1})$, we have
\begin{eqnarray*}
\Hess \gamma(q)(e,e) \leq d\sqrt{G(\gamma)}\left(\int_{0}^{\gamma}\frac{ds}{\sqrt{G(s)}} + 1\right)\ln\left(\int_{0}^{\gamma}\frac{ds}{\sqrt{G(s)}} + 1\right) + \langle \grad \gamma,\alpha(e,e)\rangle .
\end{eqnarray*}

Thus, by $(\ref{hipcurvmedia})$ it follows that
\begin{eqnarray*}
\Delta \gamma \leq B\sqrt{G(\gamma)}\left(\int_{0}^{\gamma}\frac{ds}{\sqrt{G(s)}} + 1\right)\ln\left(\int_{0}^{\gamma}\frac{ds}{\sqrt{G(s)}} + 1\right).
\end{eqnarray*}
Concluding that $(\gamma ,G)$ is an Omori-Yau pair for the laplacian in $M$.
\end{proof}

The following theorem extends the result in \cite{alias-bessa-dajczer}. We using again the following function:
\[
C_b(t)=
\begin{cases}
\sqrt{b}\cot(\sqrt{b}\,t) & \mbox{\rm if $b>0$ and $0<t<\pi/2\sqrt{b}$}\\
1/t & \mbox{\rm if $b=0$ and $t>0$}\\
\sqrt{-b}\coth(\sqrt{-b}\,t) & \mbox{\rm if $b<0$ and $t>0$}.
\end{cases}
\]

\begin{theorem}
Let $ \varphi \colon M^{m} \to L^{\ell}\times P^{n} $ be an isometric immersion where $ L^\ell $ carries an Omori-Yau pair  $ (G,\tilde{\gamma}) $ for the Hessian and $ P^n $ has a pole  $ z_0 $. If $ \pi_1 \circ \varphi $ is proper and satisfies:
\begin{itemize}
\item[i)] $\displaystyle{\lim_{x \rightarrow \infty}\frac{\rho_{P}\circ \pi_{2}(\varphi(x))}{\psi(\tilde{\gamma}\circ \pi_{1}(\varphi(x)))}} = 0$, with $ \displaystyle{\psi(t) = \log\left(\int_{0}^{t}\frac{ds}{\sqrt{G(s)}}+1\right)}$,
\item[ii)] $ \exists r > 0 $, such that $ \rho_{P}\circ \pi_{2}(\varphi(x)) \geq r $, off a compact set.
\end{itemize}
Then
\begin{eqnarray*}
\sup_{M}\vert H\vert \geq \frac{m-\ell}{m}C_{b}(\limsup_{x\to \infty}\rho_{P}\circ \pi_{2}\circ \varphi (x))
\end{eqnarray*}
where  $ \pi_1 : L^{\ell}\times P^{n} \to L^{\ell} $, $ \pi_2 : L^{\ell}\times P^{n} \to P^{n} $ are the standard projections, $ H $ is the mean curvature vector field of $ \varphi $ and $ \rho_{P} = dist_{P}(z,z_0) $.
\end{theorem}
\begin{proof}
Define $ r : L^{\ell}\times P^{n} \to \mathbb{R} $ by
\begin{eqnarray*}
r(y,z) = \rho_{P}(z)
\end{eqnarray*} 
and $ u : M^m \to \mathbb{R} $ by
\begin{eqnarray*}
u(x) = r\circ \varphi(x).
\end{eqnarray*}
We fix $ \{e_{1},...,e_{m}\} $ an orthonormal frame of $ TM $ and we write $ e_j = e_{j}^{L} + e_{j}^{P} $. In this way,
\begin{eqnarray}\label{eql1}
\Hess u(e_{i},e_{j}) = \Hess \rho_{P}(e_{i}^{P},e_{j}^{P}) + \langle \grad \rho_{P},\alpha(e_{i},e_{j})\rangle ,
\end{eqnarray}
because $ \Hess \rho_{P}(e_{i}^{P},e_{j}^{P}) = \Hess r_{L\times P}(e_{i},e_{j}) $ and $ \grad \rho_{P} = \grad_{L\times P} r $. 

For other hand,
\begin{eqnarray*}
1 = \langle e_{j},e_{j}\rangle = \vert e_{j}^{P}\vert^2 + \vert e_{j}^{L}\vert^2 .
\end{eqnarray*}
whence
\begin{eqnarray*}
m = \sum_{j=1}^{m}\left(\vert e_{j}^{P}\vert^2 + \vert e_{j}^{L}\vert^2 \right)
\end{eqnarray*}
and
\begin{eqnarray}\label{eql2}
\sum_{j=1}^{m}\vert e_{j}^{P}\vert^2 \geq (m-l).
\end{eqnarray}
From the Hessian comparison theorem applied to the manifold $ P $, we obtain
\begin{eqnarray}\label{eql3}
\Hess r(e_{j}^{P},e_{j}^{P}) \geq C_{b}(\rho_{P})(\vert e_{j}^{P}\vert^2 - \langle \grad \rho_{P},e_{j}^{P}\rangle^2).
\end{eqnarray}
Taking the trace in \eqref{eql1} and using \eqref{eql2} and \eqref{eql3}, we have
\begin{eqnarray*}
\Delta u \geq C_{b}(u)\left((m-l) - \vert \grad \rho_{P}\vert^2\right) + m\langle \grad \rho_{P},H\rangle .
\end{eqnarray*}
But $ \langle \grad \rho_{P},e_{j}^{P}\rangle = \langle \grad u,e_j \rangle $ and thus
\begin{eqnarray*}
\vert H\vert \geq \frac{m-l}{m}C_{b}(u) - \frac{1}{m}\left(C_{b}(u)\vert \grad u\vert^2 + \Delta u \right).
\end{eqnarray*}

If $ M $ is a compact manifold the result follows by computing the inequality at a point of maximum of u. Otherwise, observe that Proposition \ref{pullbackomoriyau} implies that $ M $ has an Omori-Yau pair for the Laplacian. Since 
\begin{eqnarray*}
\displaystyle{\lim_{x \rightarrow \infty}\frac{\rho_{P}\circ \pi_{2}(\varphi(x))}{\psi(\tilde{\gamma}\circ \pi_{1}(\varphi(x)))}} = 0 ,
\end{eqnarray*}
we have by the Omori-Yau maximum principle that there exists a sequence $ x_k \in M $ such that 
\begin{eqnarray*}
\vert \grad u\vert(x_k) < \frac{1}{k} \ \ \ \mbox{and} \ \ \ \Delta u(x_k) < \frac{1}{k}.
\end{eqnarray*}
Therefore,
\begin{eqnarray}
\sup_{M}\vert H\vert \geq \frac{m-l}{m}C_{b}(u(x_k)) - \frac{C_{b}(u(x_k))}{mk^2} - \frac{1}{mk}.
\end{eqnarray}
We observe that when $ b \leq 0 $, the function $ C_{b} $ is limited when $ x_k \to +\infty $. If $ b > 0 $, we have that $ 0 < t < \pi/2\sqrt{b} $ and so $ C_{b} $ also is limited by the hypothesis $ ii) $. Therefore, letting $ k \to +\infty $, we get
\begin{eqnarray*}
\sup_{M}\vert H\vert \geq \frac{m-\ell}{m}C_{b}(\lim_{k\to \infty}\rho_{P}\circ \pi_{2}\circ \varphi (x_{k}))\geq  \frac{m-\ell}{m}C_{b}(\limsup_{x\to \infty}\rho_{P}\circ \pi_{2}\circ \varphi (x)),
\end{eqnarray*}
which concludes our proof.
\end{proof}

As a last application we apply the fact that the Omori-Yau maximum principle remains valid for functions not limited that satisfying certain growth conditions. We denote $ N^{n+1} = I\times_{\rho}P^n $ the product manifold endowed with the metric of warped product, $ I \subset \mathbb{R} $ is a open interval, $ P^n $ is a complete Riemannian manifold and $ \rho : I \to \mathbb{R}_{+} $ is a smooth function. Given an isometrically immersed hypersurface  $ \psi : M^{n} \to N^{n+1} $, define $ h : M^{n} \to I $ the $ C^{\infty}(M^{n}) $ height function by setting $ h = \pi_{I}\circ f $. The result below generalizes Theorem 7 in \cite{alias-impera-rigoli}. Observe that we may have $ \sup_{M^n}h =  + \infty $.

\begin{theorem}\label{hypersufacetheorem} Let $ \psi : M^{n} \to N^{n+1} $ be an immersed hypersurface. If the Omori-Yau maximum principle holds on $ M^n $ for the Laplacian and the height function $ h $ satisfies 
$
\displaystyle{\lim_{x \rightarrow \infty}\frac{h(x)}{\varphi(\gamma(x))}} = 0 
$, where $ \varphi $ and $ \gamma $ are as in Theorem \ref{Principal Theorem}, then
\begin{equation}
\sup_{M^n}\vert H \vert \geq \inf_{M^n}\mathcal{H}(h) ,
\end{equation}
with $ H $ being the mean curvature and $ \displaystyle{\mathcal{H}(t) = \grad \ln \rho} $.
\end{theorem}


\begin{corollary} Let $ P^n $ be a complete, non-compact, Riemannian manifold whose radial sectional curvature satisfies condition 
\begin{eqnarray*}
K_{P^n}(x) \geq - B^{2}\prod_{j=1}^{k}\left[\ln^{(j)}\left(\int_{0}^{\rho(x)}\frac{ds}{\sqrt{G(s)}} + 1\right)+1\right]^{2}\left(\int_{0}^{\rho(x)}\frac{ds}{\sqrt{G(s)}} + 1\right)^{2}G(\rho(x)),
\end{eqnarray*}
where $ \rho $ is the distance function. If $ f : M^{n} \to N^{n+1} $ is a properly immersed hypersurface and the height function $ h $ satisfies the conditions imposed in Theorem \ref{hypersufacetheorem}, then
\begin{equation}
\sup_{M^n}\vert H \vert \geq \inf_{M^n}\mathcal{H}(h).
\end{equation}
\end{corollary}

\vspace{.5cm}

\noindent {\bf Acknowledgements:} We want to express our gratitude to G. Pacelli Bessa and to Newton Santos for their suggestions along the preparation on this paper. The second author would like to express his thanks to the Professor Barnabe Lima and G. Pacelli Bessa for their
 advice and illuminating discussions on Omori-Yau maximum principles.

\end{document}